\documentclass[11pt,a4paper,reqno]{amsart} %maths
\usepackage{tikz}
\usetikzlibrary{positioning,decorations.pathmorphing}

\usepackage{amssymb,graphicx}%maths
\usepackage{amssymb}
\newcommand{\koniec}{\begin{flushright}  $\Box $ \end{flushright}}
\newtheorem{theo}{Theorem}[section] 
\newtheorem{prop}[theo]{Proposition}  
\newtheorem{lemma}[theo]{Lemma}

\newtheorem{col}[theo]{Corollary}

\newcounter{mnotecount}[section]

\renewcommand{\themnotecount}{\thesection.\arabic{mnotecount}}

\newcommand{\mnote}[1]%{}%
{\protect{\stepcounter{mnotecount}}$^{\mbox{\footnotesize
$%\!\!\!\!\!\!\,
\bullet$\themnotecount}}$ \marginpar{%\color{red}%
\raggedright\tiny\em
$\!\!\!\!\!\!\,\bullet$\themnotecount: #1} }

\newcommand{\CP}{\mathbb{CP}}

\newcommand{\D}{\mathcal{D}}
\newcommand{\B}{\mathcal{B}}
\newcommand{\LL}{\mathcal{L}}
\newcommand{\CC}{\mathcal{P}}
\newcommand{\RR}{\mathcal{R}}

\newcommand{\R}{\mathbb{R}}

\def\p{\partial}
\def\be{\begin{equation}}

\def\ee{\end{equation}}

\newcommand{\pr}[2]{\langle#1, #2\rangle}
\begin{document}\date{24 September 2021}
%%%%%%%%%%%%%%%%%%%%%%%%%%%%%%%%%%%%%%%%
\vspace*{-1.0cm}
%\hfill hep-th/yymmnnn \\
%\hfill UB-ECM-PF-06-43 \\
%\hfill DMUS-MP-18-02
%\\
%\vspace{1.0cm}
\title{Variational principles for conformal geodesics}
\author{Maciej Dunajski}
\address{Department of Applied Mathematics and Theoretical Physics\\ 
University of Cambridge\\ Wilberforce Road, Cambridge CB3 0WA, UK.}
\email{m.dunajski@damtp.cam.ac.uk}

\author{Wojciech Kry\'nski}
\address{
Institute of Mathematics\\ Polish Academy of Sciences
\\Sniadeckich 8, 00-656 Warszawa, Poland
}
\email{krynski@impan.pl}
\maketitle
%%%%%%%%%%%%%%%%%%%%%%%%%%%%%%%%%%%%%%%%
%\author{Maciej Dunajski}
%\address{Department of Applied Mathematics and Theoretical Physics\\ 
%University of Cambridge\\ Wilberforce Road, Cambridge CB3 0WA, UK.}
%\email{m.dunajski@damtp.cam.ac.uk}

\begin{abstract} 
Conformal geodesics are solutions to a system of third order  equations, which makes a Lagrangian formulation problematic.
We show how enlarging the class of allowed variations leads to a variational formulation for this system with a third--order conformally invariant Lagrangian.
We also discuss the conformally invariant system of fourth order ODEs arising from this Lagrangian, and show that some of its integral curves are spirals. 
\end{abstract}   
%\maketitle
\section{Introduction}
A geodesic on a (pseudo) Riemannian manifold is uniquely specified by
a point and a tangent direction. Moreover, if two points are sufficiently close
to each other, then there exists exactly one
 length minimising (or, in Lorentzian geometry, maximising) curve between these two points, which is geodesic.  This formulation relies on the methods of
the calculus of variations, and enables a construction of normal neighbourhoods,
as well as the analysis of the Jacobi fields determining the existence of conjugate points. 

 The variational formulation has been missing in conformal geometry, where {\em conformal geodesics} arise as solutions to a system of third order ODEs: A conformal geodesic  is uniquely specified by a point, a tangent direction, and a perpendicular acceleration \cite{yano, BE, Tod}. The odd order of the underlying system of equations is not amenable to the usual methods of calculus of variation, where the resulting Euler--Lagrange equations for non--degenerate Lagrangians are of even order.

 One way around this difficulty, which we explore in this paper, is to consider a more general class of variations. As we shall see,
 this
allows to terminate the procedure of integration by parts when the integrand depends on a derivative of a variation. This argument relies on a number of technical steps: carefully controlling boundary terms, respecting conformal invariance, and making sure that the fundamental lemma of the calculus of variations can be applied to the extended class of variations.

In the next section we shall formulate the conformal geodesic equations, and summarise the notation. In \S\ref{svariations} we shall propose two ways to deal with third order equations from the variational perspective.
In \S\ref{MTs} we formulate the main result of our paper (Theorem \ref{mainth}) and compute the variation of the conformally invariant functional associated to a third order Lagrangian. While the standard variational procedure leads to a 4th order conformally invariant equation (\ref{4th_order}), looking at the extended class of variations reduces the order of the Euler--Lagrange equations to $3$, and gives conformal geodesics as extremal curves. In \S\ref{cfc} we focus on the conformally flat case, where the link between the 3rd and 4th order systems is particularly clear, and the Hamiltonian formalism can be constructed. In particular we show that logarithmic spirals arise
as solutions to the 4th order system for a particular class of initial conditions.
In \S\ref{sectiontractor} we show how the Lagrangian of Theorem \ref{mainth} can be interpreted as the `free particle' Lagrangian on the total space of the tractor bundle.
Finally in  \S\ref{dlss} we construct a degenerate Lagrangian which uses a skew--symmetric tensor.
This gives rise, via a Legendre transform, to a Poisson structure on the second--order tangent bundle.
\subsection*{Acknowledgements} 
MD has been partially supported 
by STFC grants ST/P000681/1, and  ST/T000694/1. He is also grateful to the Polish Academy of Sciences for the hospitality when some of the
 results were
obtained. WK has been supported by the Polish National Science Centre grant
 2019/34/E/ST1/00188. We thank Ian Anderson, Peter Cameron, Peter Olver, David Robinson,
 Paul Townsend and Josef Silhan for correspondence, and the anonymous
 referee for pointing out a gap in the original proof of
 Corollary \ref{maincol}.
\section{Conformal geodesic equations}
\label{section1}
A conformal class on an $n$--dimensional smooth manifold $M$ is an equivalence class of (pseudo) Riemannian metrics, where two metrics
$\hat{g}$ and $g$ are equivalent if there exists a nowhere zero function  $\Omega$ on $M$ such that
\be
\label{equivalent}
\hat{g}=\Omega^2 g.
\ee
If a metric $g\in [g]$ has been chosen, then 
${\langle }X, Y{\rangle}$ denotes the inner product of two vector fields
with respect to this metric. We also set $|X|^2\equiv \pr{X}{X}$, and use the notation
$\psi(X)$ for the $(k-1)$--form arising as a contraction of the $k$--form $\psi$ with the vector field $X$.

Let $\gamma$ be a curve of class at least $C^3$ in $M$, parametrised by $t$, and let $U$ be
a nowhere vanishing tangent vector to $\gamma$ such that $U(t)=1$. The acceleration vector of $\gamma$ is $A=\nabla_U U$, where $\nabla_U \equiv U^a\nabla_a$ is the directional derivative along $\gamma$, and $\nabla$ is the Levi--Civita connection of $g$.
The conformal geodesic equations in their conformally invariant form given by Bailey and Eastwood \cite{BE} are
\be
\label{Eeq}
E\equiv \nabla_U A-\frac{3\pr{U}{A}}{|U|^2}A+\frac{3|A|^2}{2|U|^2}U-|U|^2P^{\sharp}(U)+2P(U,U)U=0.
\ee
The Schouten tensor $P\in \Gamma (TM\otimes TM)$ of $g$ is given by 
\[
P=\frac{1}{n-2}\Big(r-\frac{1}{2(n-1)}Sg\Big),
%P_{ab}=\frac{1}{n-2}\Big(R_{ab}-\frac{1}{2(n-1)}Rg_{ab}\Big), \quad
%a, b=1, 2, \dots, n
\]
in terms of the Ricci tensor $r$ and the Ricci scalar $S$ of $g$.
The symbol $P^{\sharp}$ is an endomorphism of $TM$ defined by 
$\pr{P^{\sharp}(X)}{Y}=P(X, Y)$ for all vector fields $X, Y$.

%We will use a mixture of index--free, and the abstract index %notation, where the indices do not indicate any preferred choice of %basis or coordinates, but only indicate  the type of bundle. Thus %$U^a$ is a section of 
%$TM$, and $T_{ab}$ is a section of $T^*M\otimes T^* M$. The %exception will be our discussion of Poisson structures in \S, where %the indices to coordinates.

Changing the metric to $\hat{g}$ as in (\ref{equivalent})
results in changes to the Schouten tensor, the Levi--Civita connection
and the acceleration
\begin{eqnarray*}
\hat{P}&=&P-\nabla\Upsilon+
\Upsilon\otimes\Upsilon+|\Upsilon|^2g,\\
\hat{\nabla}_X Y&=&\nabla_X Y+\Upsilon(Y)X+\Upsilon(X)Y-\pr{X}{Y}\Upsilon^{\sharp},\\
\hat{A}&=&A-|U|^2\Upsilon^{\sharp}+2\Upsilon(U)U,
\end{eqnarray*}
where $\Upsilon\equiv \Omega^{-1}d\Omega$, and $\Upsilon^{\sharp}$ is a vector field defined by
$\Upsilon(X)=\pr{\Upsilon^{\sharp}}{X}$.
%Unlike the velocity vecotr $U$, the acceleration is not conformally invariant, and transforms according 
%\be 
%\hat{A}=A-\Upsilon  +<\Upsilon, U>U
%\ee
It is now a matter of explicit calculation to verify that the conformal geodesic equations (\ref{Eeq}) are conformally invariant.
\section{Lagrangians for third order equations}
\label{svariations}
For a non--degenerate Lagrangian, the order of the resulting Euler--Lagrange equations is equal to twice the order of the highest derivative appearing in the Lagrangian, so that in particular the Euler--Lagrange equations are of even order.   Two approaches can be taken to deal with third order systems (while they will also be applicable to systems of other odd orders, for clarity we focus on order $3$)
\begin{enumerate}
\item
To allow Lagrangians which are quadratic in the acceleration, and terminate the procedure of integration by parts at the level of 3rd order derivatives consider variations
$V$ which do not keep end points  fixed, but only satisfy $\dot{V}(t_0)=\dot{V}(t_1)=0$, where
$t_0, t_1$ are values of the parameter at end points. This enlarges the class of variations of extremal curves, 
and so reduces the number of these curves as well as the order of the resulting Euler--Lagrange equations
from $4$ to $3$.
\item Consider degenerate Lagrangians which are linear in acceleration, and necessarily involve
an anti--symmetric tensor.
\end{enumerate}
To illustrate both approaches with an elementary example, consider a curve  $\gamma$ in $\R^n$ parametrised by $t\rightarrow X(t)$, and 
aim to obtain the third order system
\be
\label{3rdo}
\dddot{X}=0
\ee
from a variational principle.  Let $\Gamma:[-1, 1]\times [t_0, t_1]\rightarrow M=\R^n$ be a
one--parameter family of curves parametrised by
$s\in[-1, 1]$,  such that 
\[
\Gamma(0, t)=\gamma(t), \quad \frac{\p \Gamma}{\p t}|_{s=0}=U, \quad \mbox{and}\quad \frac{\p \Gamma}{\p s}|_{s=0}=V 
\]
so that the variation $V$ is also a vector field along $\gamma$. 
\vskip5pt

 In the first approach we take
\[
I[X]=\frac{1}{2}\int_{t_0}^{t_1}|\ddot{X}|^2 dt
\]
so that one integration by parts gives
\[
I[X+sV]=I[X]-s\int_{t_0}^{t_1}\pr{\dddot{X}}{\dot{V}}dt+o(s).
\]
If the variation $\dot{V}$ vanishes at the end points, and is otherwise arbitrary then $\delta I=0$ iff
(\ref{3rdo}) holds.

There appear to be two immediate obstructions to generalising this approach to the conformal geodesic
system (\ref{Eeq}), where $\dot{X}=U, \ddot{X}=A+\dots$, and $\dddot{X}=\nabla_U A+\dots$, where $\dots$ denote the lower order terms involving the Christoffel symbols and curvature. Firstly, if there is an explicit $X$--dependence in the Lagrangian,
then
the undifferentiated variation $V$ appears in the integrand. Secondly
$\dot{V}=\nabla_U V+\dots$ is not conformally invariant. We shall get around both obstructions by modifying
$\nabla_U V$ to  a first order conformally invariant linear operator along $\gamma$
\be
\label{D_op}
\D(V)=\nabla_UV+|U|^{-2}(\pr{A}{V}U-\pr{U}{V}A-\pr{A}{U}V).
\ee
This operator differs from the derivative $\nabla_U$ along $\gamma$ by a linear operator which depends on the second jet of $\gamma$. It is the unique, up to the reparametrisation of $\gamma,$ conformally invariant adjustment of $\nabla_U$. 
In the Proof of Theorem \ref{mainth} we shall see that the linear operator of order zero, $D(V)-\nabla_UV$, 
has the effect of 
canceling some of the $V$ contributions in the variation of the functional, and that all these contributions can be canceled in the conformally flat case.
\vskip5pt
For the second approach, let $\Omega\in \Lambda^2(\R^n)$ be a constant two--form, and set
\[
I[X]=\int_{t_0}^{t_1}\Omega(\ddot{X}, \dot{X}) dt.
\]
If $\Gamma(s, t)=X(t)+sV(t)+o(s)$, and the variation $V$ and its derivative vanish
at the end--points, then
integrating  by parts twice give
\[
\delta I=\int_{t_0}^{t_1} 2\Omega(\dddot{X}, V)dt, \qquad\mbox{so that}\quad \Omega(\dddot{X}, \,\cdot)=0.
\]
If the dimension $n$ is even, and $\Omega$ is non--degenerate (so that $\Omega$ is a symplectic form), then
this implies (\ref{3rdo}).
\section{Main theorem}
\label{MTs}
Let  $E$ be the vector field along a smooth curve $\gamma$  defined by the equation
(\ref{Eeq}).
Define a third--order  Lagrangian $\LL$, and the corresponding
functional $I[\gamma]$ by
\be
\label{Leq}
\LL=\frac{\pr{U}{E}}{|U|^2},
\ee
and
\be
\label{functional}
I[\gamma]=\int_{t_0}^{t_1} \LL dt.
\ee
The Lagrangian $\LL$ is conformally invariant, because the expression $E$ is.
\begin{theo}\label{mainth} The first variation 
of the functional (\ref{functional}) is given by
\be
\label{main_formula}
\delta I=\int_{t_0}^{t_1} |U|^{-2}(\pr{K}{V}-\pr{E-2\LL U}{\D(V)})dt+\B(V)|^{t_1}_{t_0},
\ee
where $K$ is a vector field along $\gamma$ given, in terms of the Weyl tensor $W$, by 
%\footnote{In \cite{SZ} it is argued that
%$K=U^b({\bf \Omega}_{ab}({\bf A}, {\bf U})$, where ${\bf %\Omega}_{ab}$ is the curvature of the tractor connection.}
\be
\label{Keq}
K^e=g^{ec}({W_{bca}}^d U^a U^b A_d-2|U|^2\nabla_{[c} P_{a]b} U^aU^b),
\quad a,b, \dots=1, \dots, n
\ee
and
\be
\label{Bterm}
\B(V)=|U|^{-2}(\pr{U}{\D^2(V)}-\pr{E-2\LL U}{V}),
\ee
where $\D$ is the operator (\ref{D_op}).
\end{theo}
\noindent
{\bf Proof.} The proof is by a cumbersome calculation. We shall list the main steps, and give enough details so that the reader can verify our computations.

 The third-order term 
$|U|^{-2}\pr{U}{\nabla_U A}$ in (\ref{Leq}) can be integrated by parts to give
\[
\frac{d}{dt}\Big(\frac{\pr{U}{A}}{|U|^2}\Big)-
\Big(\frac{|A|^2}{|U|^2}-\frac{2\pr{U}{A}^2}{|U|^4}\Big), 
\]
which results in the alternative form
\be
\label{other_form}
\LL=\frac{d}{dt}\Big(\frac{\pr{U}{A}}{|U|^2}\Big)+\LL_1, \quad \mbox{where}\quad
\LL_1=\frac{1}{2}\frac{|A|^2}{|U|^2}-\frac{\pr{U}{A}^2}{|U|^4}+P(U, U).
\ee
The term $\LL_1$ coincides, up to a constant multiple, with the Lagrangian considered in \cite{BE}. Neither the resulting boundary term, nor $\LL_1$ are conformally invariant. We shall therefore focus on $\LL$, but use the variation of $\LL_1$ as an intermediate step.

 First disregard the boundary term
in (\ref{other_form}), and 
consider variations of the functional
$I_1[\gamma]=\int_{t_0}^{t_1}\LL_1 dt$. This yields
\begin{eqnarray*}
\delta I_1&=&
\int_{t_0}^{t_1}\biggl[\frac{\pr{\nabla^2_UV}{A}}{|U|^2}-\frac{\pr{A}{A}\pr{\nabla_UV}{U}}{|U|^4}+
\frac{\pr{R(V,U)U}{A}}{|U|^2}\\
&&-2\frac{\pr{\nabla_U^2V}{U}\pr{U}{A}}{|U|^4}
-2\frac{\pr{\nabla_UV}{A}\pr{U}{A}}{|U|^4}\\
&&+4\frac{\pr{U}{A}^2\pr{\nabla_UV}{U}}{|U|^6}
+(\nabla_VP)(U,U)+2P(\nabla_UV,U)\biggr]dt.
\end{eqnarray*}
The appearance of the Riemann tensor $R$ arises from varrying the metric $g$ in the inner products. We eliminate $R$ in favour of the Weyl tensor $W$ and the
Schouten tensor $P$ using the formula
$$
\begin{aligned}
\pr{W(V,U)U}{A}=&\,\,\pr{R(V,U)U}{A}-\pr{U}{U}P(V,A)-\pr{V}{A}P(U,U)\\
&+\pr{A}{U}P(V,U)+\pr{V}{U}P(A,U).
\end{aligned}
$$
We substitute
$$
2P(\nabla_UV,U)=2\nabla_U(P(V,U))-2(\nabla_UP)(V,U)-2P(V,A),
$$
and integrate the following terms, with the given coefficients, by parts:
$$
\frac{\pr{\nabla^2_UV}{A}}{|U|^2},\quad
- \frac{\pr{\nabla_UV}{A}\pr{U}{A}}{|U|^4},\quad
-\frac{1}{2}\frac{\pr{A}{A}\pr{\nabla_UV}{U}}{|U|^4},
$$
$$
-2\frac{\pr{\nabla_U^2V}{U}\pr{U}{A}}{|U|^4},\qquad 2\frac{\pr{U}{A}^2\pr{\nabla_UV}{U}}{|U|^6},\quad 
\nabla_U(P(V,U)).
$$
These terms were selected by a systematic, but somewhat tedious procedure,
which starts from the highest order term and ensures that, apart from
the inner product $\pr{K}{V}$, only terms involving $\D(V)$ appear in 
the integrand. The boundary terms arising from these integrations are combined with the 
boundary term (\ref{other_form}), which gives an expression for the variation
$\delta I$. To arrive at the statement (\ref{main_formula}) in the Theorem
we use (\ref{D_op}) to eliminate $\nabla_U V$ in favour of $\D(V)$. 
\koniec
\subsection{Conformal geodesic equations}
The boundary term  $\B$ given by (\ref{Bterm}) is conformally invariant, as both $\D$ and $E$ are. Furthermore, the vector field $K$ along $\gamma$ given by (\ref{Keq}) is conformally invariant, and consequently the integral
\be
\label{Kinn}
\mathcal{K}(V)\equiv \int_{t_0}^{t_1}|U|^{-2}\pr{K}{V}dt,
\ee
defines a conformally invariant linear operator acting on variational vector fields along a 
given curve $\gamma$. We shall exploit both $\B$ and $\mathcal{K}$ to define a class of variations needed in the following corollary
\begin{col}
  \label{maincol}
The functional $I[\gamma]$ is stationary under the class of variations
such that
\be
\label{ec1}
\B(V)|^{t_1}_{t_0}=-\mathcal{K}(V)
\ee  
if and only if  the conformal geodesic equations (\ref{Eeq}) hold.
\end{col}
\noindent
{\bf Proof.}
The proof relies on a modification of the fundamental lemma of calculus of variations, which we recall following \cite{GF}
\begin{lemma}[Fundamental Lemma of Calculus of Variations]
  If $Y:[t_0, t_1]\rightarrow\R^n$ is continuous, and such that
  \[
\int_{t_0}^{t_1}\pr{W}{Y}dt=0
\]
for all continuous $W:[t_0, t_1]\rightarrow \R^n$ with $W(t_0)=W(t_1)=0$, then $Y$ is identically $0$.
\end{lemma}
In the proof one assumes that $W$ is compactly supported with, and has
continuous derivatives up to some  specified order. For example all
components of $W$ may vanish, except the $k$th component which is given by a bump function
$\rho(t)$, where
\[
  \rho(t)=
  \begin{cases}
    (t-a)(b-t) & \text{if } t\in (a, b), \quad t_0<a<b<t_1\\
    0 &\text{otherwise.}
    \end{cases}.
\]

We now proceed to proving Corollary \ref{maincol}.
First observe that formula (\ref{main_formula}) immediately proves that $I[\gamma]$ is stationary under the class of variations such that (\ref{ec1}) holds,  provided that $\gamma$ is a solution to (\ref{Eeq}). In order to establish the converse we shall proceed by contradiction. For this assume that $\delta I[\gamma]=0$ but $E(t^*)\neq 0$ at some point $t^*\in [t_0,t_1]$ (by continuity we may assume that $t^*$ is in the interior of the interval). Our goal is to construct a variation $V$ satisfying (\ref{ec1}), and 
such that
\begin{equation}\label{eqContradiction}
\int_{t_0}^{t_1}|U|^{-2}\pr{E-2\LL U}{\D(V)}dt>0
\end{equation}
which, together with (\ref{main_formula}), and the Fundamental Lemma applied to
$W={\mathcal D}(V)$
will give a contradiction with $\delta I[\gamma]=0$.

In order to complete this task let us notice that there is a vector field $V_0$ in the kernel of the differential operator $\D$ such that $V_0(t^*)=E(t^*)-2\LL U(t^*)$. Indeed, $\D$ is a first order differential operator and one can impose arbitrary condition on the value of $V_0$ at a given point $t^*$ when solving the ODE $\D(V_0)=0$. Now, let $\rho$ be a bump function concentrated in a sufficiently small neighbourhood of $t^*$, and let $\dot \phi=\rho$. Set $\widetilde V=\phi V_0$. Then $\D(\widetilde V)=\rho V_0$ and we find  
\[
\int_{t_0}^{t_1}|U|^{-2}\pr{E-2\LL U}{\D(\widetilde V)}dt>\omega
\]
for some real number $\omega>0$. If $\widetilde V$
satisfies (\ref{ec1}), then the proof is complete. However, this condition is not satisfied in general. Therefore we shall find a correction $\hat V$ such that $V=\widetilde V+\hat V$ satisfies (\ref{ec1}) and (\ref{eqContradiction}) as well. In order to find $\hat V$ explicitly we exploit the fact that the velocity $U$
satisfies $\D(U)=0$ and $\pr{K}{U}=0$. Thus, for any function $f$ and $\hat V=fU$, we find $\mathcal{K}(\hat V)=0$ and consequently $\B(\hat V)|^{t_1}_{t_0}+\mathcal{K}(\hat V)=\left(\ddot f +f|U|^{-2}\pr{E}{U}\right)|^{t_1}_{t_0}$. Now, $f$ can be taken such that it is zero everywhere outside a small neighbourhood of $t_0$, the value of $\ddot f(t_0)$ is arbitrary large, and the values of $f$ and $\dot f$ are arbitrary small. For instance $f(t)=\frac{\kappa}{2}(t-c)^2$ for $t\in [t_0,c]$ and $f(t)=0$ otherwise satisfies these conditions for appropriately chosen constants $\kappa\in\R$ and $c>t_0$ (function $f$ can be also regularized to be of class $C^\infty$). We get that $\hat V$ can be picked such that $\B(\hat V)|^{t_1}_{t_0}+\B(\widetilde  V)|^{t_1}_{t_0}=-\mathcal{K}(\hat V)-\mathcal{K}(\widetilde V)$ on the one hand since the value of $\ddot f(t_0)$ in $\B(\hat V)|^{t_1}_{t_0}$ can take any value we want, and $\int_{t_0}^{t_1}|U|^{-2}\pr{E-2\LL U}{\D(\hat V)}dt<\frac{\omega}{2}$ on the other hand, since $\D(\hat V)=\dot f U$ and $\dot f$ can be arbitrary small. 
This completes the proof.
\koniec
\subsection{Conformal Mercator equation} Integrating 
(\ref{main_formula}) by parts once more, and neglecting the boundary term then gives
\[
\delta I=\int_{t_0}^{t_1}{\langle }-\D^*(|U|^{-2}(E-2\LL U))+|U|^{-2}K, V{\rangle}dt=0.
\]
The fundamental lemma applied to the arbitrary variation $V$ now gives the fourth order system
which we shall call the conformal Mercator equation (the terminology will be justified in the next section)
\be
\label{4th_order}
\D^*\Big(\frac{E-2\LL U}{|U|^2}\Big)-\frac{K}{|U|^2}=0,
\ee
where
\[
\D^*=-\nabla_U+|U|^{-2}({\langle }U, \cdot{\rangle}A-{\langle }A, \cdot{\rangle}U-{\langle }A, U{\rangle}\mbox{Id})
\] is the adjoint of $\D$ with respect to the $L^2$ inner product, and $K, E, \LL$ are given by
(\ref{Keq}), (\ref{Eeq}) and (\ref{Leq}) respectively.
Unlike $\D$, the operator $\D^*$ is not conformally invariant. The fourth order
system (\ref{4th_order}) is nevertheless conformally invariant
as the conformal weight $-2$ of the term $|U|^{-2}$ balances the contributions
from $\D^*$ if $g$ changes according to (\ref{equivalent}). This equation also arises from the second
order Lagrangian $\LL_1$ in (\ref{other_form}) under the assumption that the 
variation  $V$ and its derivative $\nabla_U V$ vanish at the end points $t_0$ and $t_1$. Therefore
$\LL_1$ leads to a boundary value problem for (\ref{4th_order}), where two points on the curve $\gamma$ and
two tangent vectors at these points are specified. This, in $n$ dimensions, gives $4n$ conditions which is what one would expect for a fourth order system.

If the metric $g$ is conformally flat, so that $K=0$, then any solution to the conformal geodesic equation
(\ref{Eeq}) is also a solution to (\ref{4th_order}). In the next section we shall explore this conformally flat
case in greater detail.
\section{Conformally flat case and spirals}
\label{cfc}
Assume that $M=\R^n$, the conformal class is flat, and choose a flat metric representative $g$
with vanishing Christoffel symbols. Then the Schouten tensor $P$ vanishes, and the 4th order system (\ref{4th_order})
arising, for arbitrary variations $V$ such that $V$ and its first derivatives vanish
at the end points, is
%Then the resulting Lagrangian with boundary term discarded is
%\be
%L=\frac{1}{2}\frac{|A|^2}{|U|^2}-\frac{{\langle }A, U{\rangle}^2}{|U|^4}.
%\ee
\begin{eqnarray}
\label{4th1}
\frac{dC}{dt}&=&0, \quad\mbox{where}\nonumber\\
C&=&\frac{1}{|U|^2}\Big(\dot{A}-\frac{|A|^2}{|U|^2}U-2\frac{{\langle }A, U{\rangle}}{|U|^2}A+4\frac{{\langle }A, U{\rangle}^2}{|U|^4}U
-2\frac{{\langle }\dot{A}, U{\rangle}}{|U|^2}U \Big)
\end{eqnarray}
or, eliminating the ${\langle }\dot{A}, U{\rangle}$ term, 
\be
\label{4th_2}
\dot{A}+\frac{|A|^2}{|U|^2}U-\frac{2{\langle }A, U{\rangle}}{|U|^2}A+2{\langle }C, U{\rangle}U-|U|^2C=0, \quad C=\mbox{const}.
\ee
Picard's existence and uniqueness theorem implies that a solution curve
to (\ref{4th1}) is determined by specifying initial conditions
$X(0), U(0), A(0), \dot{A}(0)$. Specifying these values also determines
$C$ in (\ref{4th1}), and conversely specifying
$X(0), U(0), A(0), C$ determines $\dot{A}(0)$ by (\ref{4th_2}). There is an advantage  in using
$C$ instead of $\dot{A}(0)$ in the initial conditions, as $C$ stays constant
along the integral curves of  (\ref{4th1}). This can be used
to verify directly that all conformal geodesics are integral curves of (\ref{4th1})
for special initial conditions. Indeed, setting
\be
\label{C_formula}
C=\Big(\frac{1}{2}\frac{|A|^2}{|U|^4}-2\frac{{\langle }A, U{\rangle}^2}{|U|^6}\Big)U+\frac{{\langle }A, U{\rangle}}{|U|^4}A,
\ee
and substituting this into (\ref{4th_2}) gives the conformal geodesic equations for the flat metric
\be
\label{flateq}
\dot{A}-\frac{3{\langle }A,  U{\rangle}}{|U|^2}A+\frac{3|A|^2}{2|U|^2}U=0.
\ee
We also verify that $\dot{C}=0$, as a consequence of (\ref{flateq}). Thus 
(\ref{C_formula}) gives a first integral
of the conformal geodesic equations. The 
general solutions to these equations
are projectively parametrised circles
\be
\label{circles1}
t\rightarrow X(t)=X_0+\frac{tU_0+t^2 A_0}{1+t^2|A_0|^2},
\ee
where $U_0$ is a constant unit vector, and ${\langle }U_0, A_0{\rangle}=0$.
Therefore (\ref{C_formula}) evaluated at $t=0$ defines a submanifold
in the space of initial conditions which singles out conformal geodesics as integral curves. 

 For generic initial conditions the integral curves of (\ref{4th1})
are not conformal geodesics.  For example, choosing arbitrary values
of $X(0), U(0), A(0)$, and setting $C=0$ reduces
(\ref{4th1}) to
\be
\label{seq11}
\dot{A}-\frac{2{\langle }A, U{\rangle}}{|U|^2}A+\frac{|A|^2}{|U|^2}U =0.
\ee
The general solution of this system 
is\footnote{To find this solution set $u=U/|U|$ so that $U=\dot{h}(t) u$, where $\dot{h}(t)=|U|$, and $s=h(t)$ is the arc-length. Substituting this into (\ref{4th_2}) eventually gives
\be
\label{44}
\ddot{u}+\Big(H^2-3\dot{H}+m\Big)u-\dot{h}C=0, \quad\mbox{and}\quad
{\langle }C, u{\rangle}=-\dot{H}/\dot{h}
\ee
where $H(t)\equiv\ddot{h}/\dot{h}$, and $m$ is a constant. The function $h(t)$ is constrained by
a scalar nonlinear  ODE which arises by dotting (\ref{44}) with $C$. For any given solution 
of this ODE, the condition (\ref{44}) becomes a linear equation for $u$, and its solution
needs to be integrated to recover the curves $t\rightarrow X(t)$. If $C=0$ then $H$ is a constant, and
an affine transformation of $t$ can be used to set $\dot{h}=e^t$. The general solution
of (\ref{44}) then gives the spirals (\ref{spirals}).}
\be
\label{spirals}
t\rightarrow X(t)=e^t\cos{(ct)}\;P_0 +e^t\sin{(ct)}\;Q_0+R_0,
\ee
where $P_0, Q_0, R_0$ are constant vectors such that ${\langle }P_0, Q_0{\rangle}=0$ and
$|P_0|=|Q_0|$. The curves (\ref{spirals}) are logarithmic spirals
in the plane spanned by $(P_0, Q_0)$ which spiral towards $R_0$ as
$t\rightarrow -\infty$. If $(r, \theta)$ are plane polar coordinates in the plane spanned by
$P_0, Q_0$ and centered at $R_0$, then the unparametrised form of the spirals are
$
r=|P_0|e^{\theta/c}.
$
This contrasts with the behaviour of conformal geodesics, where the spirals
are conjectured not to arise even in curved cases \cite{Tod, CDT}.
The conformal invariance of the fourth order system (\ref{4th_order}) ensures that the inverse images of the
logarithmic spirals under the stereographic projection from $S^{n}$ to $\R^n$ are solutions
to (\ref{4th_order}) on the round sphere. These curves are the loxodromes. They cut all meridians at a fixed angle, and correspond to straight lines on the Mercator map - this justifies our terminology.

For general initial conditions the integral curves (\ref{4th1}) are, unlike
(\ref{circles1})  and  (\ref{spirals}), no-longer
planar, and their Serret--Frenet torsion is given in
terms of different initial jerk  $\dot{A}(0)$ (see \cite{DG} for other occurrences of equations
involving a change of acceleration in physics).
\begin{center}
\label{Fig1}
\includegraphics[width=7cm,height=7cm,angle=0]{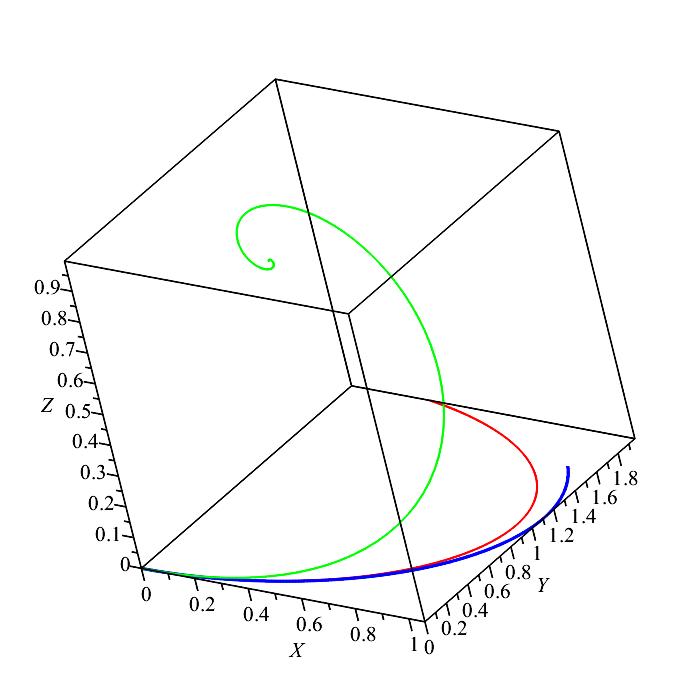}
\begin{center}
{\em {\bf Figure 1.} Integral curves of equation (\ref{4th1}) with the same initial values
of $X=(0, 0, 0), \dot{X}=(1, 0, 0), \ddot{X}=(0.1, 1, 0)$
but different $\dddot{X}$: Conformal geodesic (red), logarithmic spiral (blue), and a generic integral curve
with non--zero torsion (green).}
\end{center}
\end{center}
 In Figure 1 we plot integral curves
of (\ref{4th1}) in $\R^3$ corresponding to 3 sets of initial conditions. 
Each set shares the same values of 
$X(0), U(0), A(0)$ (so the three integral curves
are tangent to the second order at $X(0)$) but has different $\dot{A}(0)$.
If $\dot{A}(0)$ is determined by
(\ref{flateq}) in terms of the remaining initial conditions, then the
(red) integral curves are circles. If $\dot{A}(0)$ is determined
by $C=0$, or equivalently by (\ref{seq11}) then the (blue) curves are spirals.
Finally if $C=(0, 0, 1)$ then the numerical solution of (\ref{4th1})
yields the non planar (green) integral curves. The general solution to the conformal Mercator equation 
(\ref{4th1}) is given by the special conformal transformation of the logarithmic spiral
(\ref{spirals}):
\be
\label{final_s}
t\rightarrow Y(t)=\frac{X(t)-|X(t)|^2B}{1-2{\langle }X(t), B{\rangle}+|B|^2 |X(t)|^2},
\ee
where $X(t)$  is given by (\ref{spirals}), and $B$ is a constant vector\footnote{After this paper appeared on the arXiv, Josef Silhan has pointed out to us that the forthcoming work of his and Vojtech Zadnik characterise these curves by their constant conformal curvature, and vanishing conformal torsion.}.
\subsection{Hamiltonian formalism}
The fourth order system (\ref{4th1}) arises from a Hamiltonian, and we shall construct the Hamiltonian formulation using the Ostrogradsky approach to higher derivative Lagrangians.
Neglecting the boundary term in (\ref{Leq}) gives the
Lagrangian
\[
\LL=\frac{1}{2}\frac{|\dot{U}|^2}{|U|^2}-\frac{{\langle }U, \dot{U}{\rangle}}{|U|^4}-{\langle }\lambda, U-\dot{X}{\rangle},
\]
where the curve $\gamma$ is parametrised by $t\rightarrow X(t)$ with $X\in M=\R^n$, and 
$\lambda\in \R^n$ is a Lagrange multiplier. Define the momenta 
$\CC$ and $\RR$ conjugate to
$X$ and $U$ by
\[
\CC\equiv\frac{\p \LL}{\p \dot{X}}=\lambda, \quad \RR\equiv \frac{\p \LL}{\p \dot{U}}=\frac{\dot{U}}{|U|^2}-2\frac{{\langle }U, \dot{U}{\rangle}}{|U|^4}U,
\]
so that $\CC-\lambda=0$ is the set of constraints, and we can eliminate $\dot{U}$ by
\[
\dot{U}=|U|^2(\RR-2|U|^{-2}{\langle }U, \RR{\rangle}U).
\]
In this formula, and below, we abuse notation and use the flat metric to identify the tangent space to $\R^n$ with its dual, and so regard
$\CC$ and $\RR$ as both vectors and co-vectors depending on the context.
The Hamiltonian is given by the Legendre transform
\begin{eqnarray}
H&=&{\langle }\CC, \dot{X}{\rangle}+{\langle }\RR, \dot{U}{\rangle}-\LL\nonumber\\
&=&\frac{1}{2}|U|^2|\RR|^2-{\langle }U, \RR{\rangle}^2+{\langle }\CC, U{\rangle}.
\end{eqnarray}
Using the canonical commutation relations
\[
\{X^a, \CC_b\}={\delta^a}_b, \quad \{U^a, \RR_b\}={\delta^a}_b
\]
gives 
\begin{eqnarray*}
\dot{X}^a&=&U^a, \quad \dot{U}^a=|U|^2\RR^a-2{\langle }U, \RR{\rangle}U^a, \\
\dot{\RR}^a&=&\{\RR^a, H\}=-|\RR|^2U^a+2{\langle }U, \RR{\rangle}\RR^a-\CC^a, \quad \dot{\CC}^a=0
\end{eqnarray*}
Turning this system of 1st order ODEs to a 4th order system yields
(\ref{4th_2}).
%%%%%%%%%%%%%%%%%%%%%%%%%%%%%%%%%%%%%%%%%%%%%%%%%%%%%%%%
\section{Free particle on the tractor bundle}
\label{sectiontractor}
In the tractor approach of \cite{BEG}, the condition for a curve $\gamma$ to be a conformal geodesic is shown to be equivalent to 
the condition that the acceleration tractor is constant along this curve, 
and that its
tractor norm is zero. This, at least formally, has a simple mechanical interpretation of a free particle on the total space of the tractor bundle, whose position is given by the velocity tractor ${\bf U}$. The natural kinetic Lagrangian given by the squared tractor norm of the acceleration tractor is equal to 
(\ref{Leq}).

The details are as follows.
The tractor bundle is a rank $(n+2)$
vector bundle ${\bf T}={\mathcal E}[1]\oplus TM\otimes{\mathcal E}[-1]\oplus{\mathcal E}[-1]$, where
${\mathcal E}[k]$ denotes a line bundle over $M$ of conformal densities of weight $k$. 
Under the conformal rescallings (\ref{equivalent}) a section ${\bf X}=(\sigma, \mu^a, \rho)$ of ${\bf T}$ transforms according to
\[
\left(\begin{array}{c}
\hat{\sigma}\\ 
\hat{\mu}^a\\
\hat{\rho}
\end{array}\right)= \left(\begin{array}{c}
{\sigma}\\ 
{\mu}^a+\Upsilon^a \sigma\\
{\rho}-\Upsilon_a\mu^a-\frac{1}{2}|\Upsilon|^2\sigma
\end{array}\right).
\]
The tractor bundle comes equipped with a conformally invariant connection
\be
\label{tractor_con}
{\quad{{D}}_a \left(\begin{array}{c}
\sigma\\ 
\mu^b\\
\rho
\end{array} \right)= 
\left(\begin{array}{c} \nabla_a\sigma-\mu_a \\ 
\nabla_a\mu^b+P_{a}^b\sigma+\delta_a^b\rho\\
\nabla_a\rho-P_{ab}\mu^b
\end{array} \right),}
\ee
and a metric on the fibres of ${\bf T}$ defined by the norm
\[
{\langle }{\bf X}, {\bf X}{\rangle}_{\bf T}\equiv |\mu|^2+2\sigma\rho,
\]
and preserved by the tractor connection.

With the velocity tractor ${\bf U}$, and the acceleration tractor ${\bf A}$ given by
\[
{\bf U}=
\left(\begin{array}{c}
0\\ 
|U|^{-1}U\\
-|U|^{-3}{\langle }U, A{\rangle}
\end{array}\right), \quad {\bf A}=U^aD_a {\bf U}
\]
a curve is a projectively parametrised conformal geodesic
if
\be
\label{tractor}
{\langle }{\bf A}, {\bf A}{\rangle}_{\bf T}=0, \quad \frac{d {\bf A}}{d t}=0,
\ee
where now $d/dt$ is the directional tractor derivative. Given that ${\bf A}=\dot{\bf U}$, the second equation in (\ref{tractor}) is the
Euler--Lagrange equation $\ddot{\bf U}=0$ for the `free particle', with position ${\bf U}$. This equation arises from a Lagrangian $1/2{\langle }\dot{\bf U}, \dot{\bf U}{\rangle}_{\bf T}$, and
the first equation in (\ref{tractor}) states that this Lagrangian,
or equivalently the kinetic energy taken with respect to the tractor metric, is zero. It can be verified by explicit calculation that this kinetic energy is equal to the 3rd order Lagrangian (\ref{Leq}), which, at least at this formal level, therefore appears to be  natural.
\section{Degenerate Lagrangian}
\label{dlss}
The  second approach to the Lagrangian formulation alluded to in \S\ref{svariations} is to allow general variations, but consider
a degenerate Lagrangian linear in the acceleration. As explained in \S{\ref{svariations}} such Lagrangian must involve a preferred
anti-symmetric tensor, which in general is absent in conformal geometry. We shall therefore restrict to 
structures which admit a K\"ahler metric in the conformal class. In this case it is advantageous
to use a formulation of the conformal geodesic equations which is due to Yano \cite{yano},  and Tod \cite{Tod}.

This formulation
is equivalent to (\ref{Eeq}) after a change of parametrisation.
Decomposing the LHS $E$ (\ref{Eeq}) into parts orthogonal and parallel to $U$
gives $E\wedge U$, and ${\langle }E, U{\rangle}$. The second term can be made zero by
reparametrizing, and using $s=s(t)$ as a parameter along $\gamma$. An explicit calculation then verifies that vanishing of the first term $E\wedge U$ is
invariant. If one takes  $s$ to be the arc-length so that $|d\gamma/ds|^2=|U|^2=1$, and  $A=\nabla_U U$ is orthogonal to $U$, then the conformal geodesic equations become
\be
\label{pauls}
\nabla_U A=-(|A|^2+P(U, U))U+P^{\sharp}(U).
\ee
See \cite{BE, Tod, EZ20} for details.
\begin{lemma}
K\"ahler-magnetic geodesics on a K\"ahler--Einstein manifold are also conformal geodesics.
\end{lemma}
\noindent
{\bf Proof}.
Assuming that the metric $g$ is Einstein reduces (\ref{pauls}) to
\[
\nabla_U A=-|A|^2U,
\]
and the norm $|A|^2$ of the acceleration is a first integral.
If $g$ is in addition K\"ahler with the complex structure $J$, then the 
K\"ahler magnetic geodesics (i. e. the geodesics of charged particles moving  in a magnetic field given by the
K\"ahler form) are also conformal geodesics\footnote{The converse statement is not true: there are more unparametrised 
conformal geodesics (a $(3n-3)$--dimensional family)
than unparametrised K\"ahler--magnetic geodesics (a $(2n-1)$--dimensional
family, if one regards the charge as one of the parameters)
through each point. For example \cite{jap0} on $\CP^2$ 
magnetic--K\"ahler
geodesics lift to circles on $S^5$, and all conformal geodesics
lift to helices of order 2, 3 and 5. It is however the case that on a hyper--K\"ahler four--manifold
every conformal geodesic is a K\"ahler--magnetic geodesic for some choice of  K\"ahler structure 
\cite{DT_circles}.}.  Indeed, the equation
\be
\label{magneticge}
\nabla_U U=eJ(U),
\ee
where the constant $e$ is the electric charge, implies
\begin{eqnarray*}
\nabla_U A&=& \nabla_U\nabla_U U=e J(\nabla_U U)=e^2 J^2 (U)\\
&=&-e^2 U
\end{eqnarray*}
where we used $\nabla  J=0$. This is the conformal geodesic equation with $|A|=e$. 
\koniec
The K\"ahler magnetic geodesics arise from a Lagrangian
\[
L=\frac{1}{2}|U|^2+\phi(U),
\]
where the one--form $\phi$ is the magnetic potential, i. e. $d\phi=\Omega$, and $\phi(U)$ is the contraction of the vector field $U$ with $\phi$. In the corresponding Hamiltonian formalism the K\"ahler--magnetic geodesics are integral curves
of a Hamiltonian vector field on the phase--space.

A choice for the phase space in the conformal geodesic context is the second order tangent bundle $T^2 M$. It is the union of all
second order tangent spaces $T^2_xM$ - the space of equivalence
classes (also called 2-jets $j^2(\gamma)$) of curves $\gamma:\R\rightarrow M$ which agree at $x\in M$ up to and including the second derivatives.
The space $T^2M=\cup_{x\in M} T^2_x M$ is a bundle, but not, in general,
a vector bundle over $M$. 
%If $({\mathcal U}, x^i), i=1, \dots, n$ is a coordinate neighbourhood of $x\in M$ with
%$\gamma(0)=x$, then\md{skip index}
%\[
%x^i =\gamma^i(0), \quad u^i=\dot{\gamma}^i(0), \quad w^i=\ddot{\gamma}^i(0)
%\]
%are local coordinates on pre-immage of ${\mathcal U}$ in $T^2M$. 
%The coordinates $u^i$ are components of the velocity vector transform as a rank-one tensor under diffeomorphisms of $M$, but the coordinates $w^i$ transform inhomogeneously.
In the presence of a linear connection there exists a canonical splitting
\cite{yano2}
\[
S:T^2 M\rightarrow TM\oplus TM, \quad
S(j^2(\gamma))\rightarrow(\dot{\gamma}(0), (\nabla_{\dot{\gamma}(0)}\dot{\gamma})(0)).
\]
%which in local coordinates $(x, u, w)$ takes the form
%\[
%S(x^i, u^i, w^i)=(x^i, u^i, a^i\equiv w^i+{\Gamma^i}_{jk}u^ju^k).
%\]
This splitting equips $T^2M$ with the structure of rank--$2n$ vector bundle over $M$, and allows a definition of the acceleration of the curve at $\gamma(0)$ as $A=\nabla_{\dot{\gamma}(0)}\dot{\gamma}(0)$.

We would like to construct a Lagrangian on $T^2M$
involving $\phi$, as well as $\Omega$ which gives rise to all conformal geodesics on
a K\"ahler--Einstein manifold $(M, g, \Omega)$, but even this appears to be out of reach.
 The following construction works in  the flat case with $M=\R^n$, where $n$ is even, and $\Omega$ is a (chosen) K\"ahler form. 

Consider a second--order Lagrangian
\be
\label{LM1}
%L=w^2\phi_i \dot{x}^i+\frac{1}{2}\Omega_{ij}\ddot{x}^i\dot{x}^j, 
L=w^2\phi(\dot{X})+\frac{1}{2}
\Omega(\ddot{X}, \dot{X})
\ee
where $d\phi=\Omega$, or $\p_a\phi_b-\p_b\phi_a=\Omega_{ab}$, and
$w$ is a constant. 
The Euler--Lagrange equations are
\begin{eqnarray*}
0&=&
\frac{\p L}{\p x^a}-\frac{d}{ds}\frac{\p L}{\p \dot{x}^a}+
\frac{d^2}{d s^2}\frac{\p L}{\p\ddot{x}^a}\\
&=&w^2\dot{x}^b\p_a\phi_b-\frac{d}{ds}(w^2\phi_a+\frac{1}{2}\Omega_{ba}\ddot{x}^b)+\frac{d^2}{ds^2}\Big(\frac{1}{2}\Omega_{ba}\dot{x}^b\Big)\\
&=&w^2(\p_a\phi_b-\p_b\phi_a)\dot{x}^b+\Omega_{ab}\dddot{x}^b
\end{eqnarray*}
which gives, with $x^a, a=1, \dots, n$ denoting the components of $X\in \R^n$,
\be
\label{eq1}
%\dddot{x}^i=-w^2\dot{x}^i, 
\dddot{X}=-w^2 \dot{X}
\ee
as $\Omega_{ab}$ is invertible.
\begin{prop}
The Poisson structure on $T^2(\R^n)$ with $n$ even, and  coordinates $(x^a, U^a, A^a)$
induced by the Lagrangian (\ref{LM1}) is
\be
\label{poisson1}
\{x^a, A^b\}=-\Omega^{ab}, \quad 
\{U^a, U^b\}=\Omega^{ab}, \quad \{A^a, A^b\}=w^2\Omega^{ab}, 
\ee
where $w$ is a non--zero constant, and all other brackets vanish. 
Then (\ref{eq1}) are Hamilton's equations with the Hamiltonian
\be
\label{Ham_au}
H=\Omega(A, U).
\ee
\end{prop}
\noindent
{\bf Proof.} 
%Symplectic structure on the $3n$ (with even $n$)--dimensional phase space\md{Darboux coords?}
%\[
%{\bf \Omega}=\frac{1}{2}\Omega_{ij}(-w^2dx^i\wedge dx^j-2dx^i\wedge dA^j
%+dU^i\wedge dU^j).
%\]
In order to construct the canonical formalism, rewrite (\ref{LM1}) as a first order Lagrangian with constraints
\[
L=w^2\phi(U)+\frac{1}{2}\Omega(\dot{U}, U)-{\langle }\lambda, U-\dot{X}{\rangle}
%L=w^2\phi_i %U^i+\frac{1}{2}\Omega_{ji}\dot{U}^jU^i-\lambda_i(U^i-\dot{x}^i)
\]
with $n$ Lagrange multipliers $\lambda\in\R^n$. The conjugate momenta
\[
\frac{\p L}{\p \dot{x}}= \CC, 
\quad 
\frac{\p L}{\p \dot{U}}= 
\RR
\quad
\frac{\p L}{\p \dot{\lambda}}= {\mathcal S}.
\]
This gives rise to $3n$ constraints
\be
\label{constraints1}
\CC_a-\lambda_a=0, \quad \psi_a\equiv \RR_a-\frac{1}{2}\Omega_{ab}U^b=0, \quad {\mathcal S}^a=0.
\ee
Impose these constraints in the $6n$--dimensional phase-space
with coordinates $(x, \CC, U, \RR, \lambda, {\mathcal S})$ with the Poisson structure
\[
 \{x^a, \CC_b\}={\delta^a}_b, \quad
\{U^a, \RR_b\}={\delta^a}_b, \quad \{\lambda_a, {\mathcal S}^b\}={\delta^b}_a
\]
(and all other brackets vanishing),
and compute the Dirac brackets \cite{dirac} on the $3n$--dimensional reduced phase space
$TM\oplus T^*M$ with coordinates $(x, U, \CC)$. 

 Modifying $\{U^a, U^b\}=0$ to the Dirac bracket gives
\[
\{U^a, U^b\}^*\equiv\{U^a, U^b\}-\{U^a, \psi_c\}(C^{-1})^{cd}\{\psi_d, U^b
\},
\]
where $
C_{cd}\equiv\{\psi_c, \psi_d\}=\Omega_{cd}$ and
$
(C^{-1})^{dc}=\Omega^{cd}$, where $\Omega_{ab}\Omega^{ac}={\delta_b}^c$.
We also find $\{U^a, \psi_c\}={\delta^a}_c$, and finally (dropping  $*$)
\be
\label{poisson2}
\{x^a, \CC_b\}={\delta^a}_b, \quad \{U^a, U^b\}=\Omega^{ab}.
\ee
The Hamiltonian is now given by the Legendre transform (see \cite{masterov} for other possible choices of phase spaces)
\begin{eqnarray}
\label{hamiltonian}
H&=&={\langle }\dot{X},\CC{\rangle}+{\langle }\dot{U}, \RR{\rangle} +{\langle }\dot{\lambda}, {\mathcal S}{\rangle}-L\nonumber\\
&=&{\langle }\CC-w^2\phi, U{\rangle}
%H&=& \dot{x}^i\CC_i+\dot{u}^i\RR_i+ \dot{\lambda}_is_i-L=(\lambda_i-w^2\phi_i)u^i\nonumber\\
%&=&(\CC_i-w^2\phi_i)U^i,
\end{eqnarray}
where the last expression holds on the surface of constraints.
Hamilton's equations  are equivalent to (\ref{eq1}). 
%\begin{eqnarray*}
%\dot{x}^i&=&\{x^i, H\}=u^i\\
%\dot{u}^i&=&\{u^i, H\}=(p_j-w^2\phi_j)\Omega^{ij}\\
%\dot{p}_i&=&\{p_i, H\}=-w^2u^j\{p_i, \phi_j\}=w^2u^j\p_i\phi_j.
%\end{eqnarray*}
%We verify that they are equivalent to (\ref{eq1}):
%\begin{eqnarray*}
%\dddot{x}^i&=&(w^2u^k\p_j\phi_k-w^2\p_k\phi_j u^k)\Omega^{ij}\\
%&=&w^2\Omega_{jk}\Omega^{ij}u^k=-w^2x^i.
%\end{eqnarray*}
If we use the equation $\dot{U}=\{U, H\}$ to define 
$A^a=(\CC_b-w^2\phi_b)\Omega^{ab}$, and 
instead use $TM\oplus TM$ as the phase-space with coordinates
$(x^a, U^a, A^a)$, then eliminating $\CC$ by
\[
\CC_b=\Omega_{ab}A^a+w^2\phi_b
\]
yields the Hamiltonian
(\ref{Ham_au}).
The Poisson brackets (\ref{poisson2}) yield the Poisson structure (\ref{poisson1}). 
Hamilton's equations
\[
\dot{x}^a=\{x^a, H\}=U^a, \quad \dot{U}^a=\{U^a, H\}=A^a, 
\quad \dot{A}^a=\{A^a, H\}=\Omega_{bc}U^c\{A^a, A^b\}=-w^2U^a,
\]
are equivalent to (\ref{eq1}).

\koniec
The Poisson structure (\ref{poisson1}) does not generalise to curved spaces, as the Jacobi identity is obstructed by the Riemann curvature of $g$. It is nevertheless possible to make contact with the first integrals of the conformal geodesic equations, and the conformal Killing--Yano tensors
under the additional assumption that $w^2=|A|^2$.
%Symplectic structure on the $3n$ (with even $n$)--dimensional phase space\md{Darboux coords?}
%\[
%{\bf \Omega}=\frac{1}{2}\Omega_{ij}(-w^2dx^i\wedge dx^j-2dx^i\wedge dA^j
%+dU^i\wedge dU^j).
%\]
%Conventions
%\[
%X_H\hook {\bf\Omega}=dH, \quad {\bf\Omega}(X_F, X_G)=\{F, G\}=X_G(F)=X_G\hook dF.
%\]
In this case the Hamiltonian vector field corresponding to (\ref{Ham_au})
is given by
\[
X_H=U^a\frac{\p}{\p x^a}+A^a\frac{\p}{\p U^a}-|A|^2U^a\frac{\p}{\p A^a}
\]
Consider a function $Q:T^2 M\rightarrow \R$ of the form
\[
Q=Y(A, U)+W(U),
\]
where $Y\in\Lambda^1(M)\otimes\Lambda^1(M)$ and $W\in\Lambda^1(M)$ are differential forms
on $M=\R^n$ which are not necessarily constant.
The function $H$ Poisson-commutes with $Q$ iff the conformal-Killing-Yano (CKY) equation 
\be
\label{CKY}
\nabla_aY_{bc}=\nabla_{[a}Y_{bc]}+2g_{a[b}W_{c]}
\ee
holds. Indeed
\[
\{Q, H\}=X_H(Q)=(\p_a Y_{bc}+2g_{a[c}W_{b]})U^aU^cA^b+U^aU^b\p_a W_b+Y_{bc}(A^bA^c-|A|^2U^bU^c)
\]
so $\{Q, H\}=0$ iff 
\[
Y_{(bc)}=0, \quad \p_{(b} W_{c)}=0, \quad\p_a Y_{bc}= \p_{[a}Y_{bc]}+ 2g_{a[b}W_{c]}.
\]
Therefore $Q$ is is constant along the conformal geodesics, in agreement with the results of \cite{Tod} and 
\cite{rodar}.
\section{Conclusions}
Conformal geodesics are examples of distinguished curves in parabolic
geometries \cite{cap, DZ,  EZ20, SZ}. They also arise naturally in 
General Relativity
as a tool in studying the proprieties of conformal infinity \cite{Friedrich,
LT, mariem}. Despite their importance, there are few explicit examples
known \cite{Tod, DT_circles}, and the underlying mathematical theory is not
nearly as well  developed as that of geodesics. In this paper we have
proposed a variational formulation of the conformal geodesic
equations. We hope that this will shed light on the integrability
properties of these equations \cite{rodar}, as well as the global problems
such as trapping and spirals in conformal geometry \cite{CDT}.

\end{document}